# Bipolar Neutrosophic Sets And Their Application Based On Multi-Criteria Decision Making Problems

*Irfan Deli[1], Mumtaz Ali[2] and Florentin Smarandache[3]*

[1] Muallim Rıfat Faculty of Education, Kilis 7 Aralık University, 79000 Kilis, Turkey, irfandeli@kilis.edu.tr

[2] Department of Mathematics, Quaid-e-Azam University Islamabad. mumtazali7288@gmail.com

[3] University of New Mexico, 705 Gurley Ave., Gallup, New Mexico 87301, USA. fsmarandache@gmail.com

## ABSTRACT

In this paper, we introduce concept of bipolar neutrosophic set and its some operations. Also, we propose score, certainty and accuracy functions to compare the bipolar neutrosophic sets. Then, we develop the bipolar neutrosophic weighted average operator ($A_w$) and bipolar neutrosophic weighted geometric operator ($G_w$) to aggregate the bipolar neutrosophic information. Furthermore, based on the $A_w$ and $G_w$ operators and the score, certainty and accuracy functions, we develop a bipolar neutrosophic multiple criteria decision-making approach, in which the evaluation values of alternatives on the attributes take the form of bipolar neutrosophic numbers to select the most desirable one(s). Finally, a numerical example of the method was given to demonstrate the application and effectiveness of the developed method.

**Keywords:** Neutrosophic set, bipolar neutrosophic set, average operator, geometric operator, score, certainty and accuracy functions, multi-criteria decision making.

## 1. Introduction

To handle with imprecision and uncertainty, concept of fuzzy sets and intuitionistic fuzzy sets originally introduced by Zadeh [26] and Atanassov [1], respectively. Then, Smarandache [17] proposed concept of neutrosophic set which is generalization of fuzzy set theory and intuitionistic fuzzy sets. These sets models have been studied by many authors; on application [4-6,10-12,15,16], theory [18-20,21-25,27,28], and so on.

Bosc and Pivert [2] said that "Bipolarity refers to the propensity of the human mind to reason and make decisions on the basis of positive and negative affects. Positive information states what is possible, satisfactory, permitted, desired, or considered as being acceptable. On the other hand, negative statements express what is impossible, rejected, or forbidden. Negative preferences correspond to constraints, since they specify which values or objects have to be rejected (i.e., those that do not satisfy the constraints), while positive preferences correspond to wishes, as they specify which objects are more desirable than others (i.e., satisfy user wishes) without rejecting those that do not meet the wishes." Therefore, Lee [8,9] introduced the concept of bipolar fuzzy sets which is an generalization of the fuzzy sets. Recently, bipolar fuzzy models have been studied by many authors on algebraic structures such as; Chen et. al. [3] studied of $m$-polar fuzzy set and illustrates how many concepts have been defined based on bipolar fuzzy sets. Then, they examined many results which are related to these concepts can be generalized to the case of $m$-polar fuzzy sets. They also proposed numerical examples to show how to apply $m$-polar fuzzy sets in real world problems. Bosc and Pivert [2] introduced a study is called bipolar fuzzy relations where each tuple is associated with a pair of satisfaction degrees. Manemaran and Chellappa [14] gave some applications of bipolar fuzzy sets in groups are called the bipolar fuzzy groups, fuzzy d-ideals of groups under (T-S) norm. They investigate some related properties of the groups and introduced relations between a bipolar fuzzy group and bipolar fuzzy d-ideals. Majumder [13] proposed bipolar valued fuzzy subsemigroup, bipolar valued fuzzy bi-ideal, bipolar valued fuzzy (1,2)- ideal and bipolar valued fuzzy ideal. Zhou and Li [29] introduced a new framework of bipolar fuzzy subsemirings and bipolar fuzzy ideals which is a generalization of fuzzy subsemirings and bipolar fuzzy ideals in semirings and and bipolar fuzzy ideals, respectively, and related properties are examined by the authors.

In this paper, we introduced the concept of bipolar neutrosophic sets which is an extension of the fuzzy sets, bipolar fuzzy sets, intuitionistic fuzzy sets and neutrosophic sets. Also, we give some operations and operators on the bipolar neutrosophic sets. The operations and operators generalizes the operations and operators of fuzzy sets, bipolar fuzzy sets, intuitionistic fuzzy sets and neutrosophic sets which has been previously proposed. Therefore, in section 2, we introduce concept of bipolar neutrosophic set and its some operations including the score, certainty and accuracy functions to compare the bipolar neutrosophic sets. In the same section, we also develop the bipolar neutrosophic weighted average operator ($A_w$) and bipolar neutrosophic weighted geometric operator($G_w$) operator to aggregate the bipolar neutrosophic information. In section 3, based on the $A_w$ and $G_w$ operators and the score, certainty and accuracy functions, we develop a bipolar neutrosophic multiple criteria decision-making approach, in which the evaluation values of alternatives on the attributes take the form of bipolar neutrosophic numbers to select the most desirable one(s) and give a numerical example of the to demonstrate the application and effectiveness of the developed method. In last section, we conclude the paper.

## 2. Preliminaries

In the subsection, we give some concepts related to neutrosophic sets and bipolar sets.

**Definition 2.1.** [17] Let X be a universe of discourse. Then a neutrosophic set is defined as:

$$A = \{\langle x, F_A(x), T_A(x), I_A(x)\rangle : x \in X\},$$

which is characterized by a truth-membership function $T_A: X \to ]0^-, 1^+[$, an indeterminacy-membership function $I_A: X \to ]0^-, 1^+[$ and a falsity-membership function $F_A: X \to ]0^-, 1^+[$.

There is not restriction on the sum of $T_A(x)$, $I_A(x)$ and $F_A(x)$, so $0^- \leq \sup T_A(x) \leq \sup I_A(x) \leq \sup F_A(x) \leq 3^+$.

For application in real scientific and engineering areas, Wang et al.[18] proposed the concept of an single valued neutrosophic set as follows;

**Definition 2.1.** [18] Let X be a universe of discourse. Then a single valued neutrosophic set is defined as:

$$A_{NS} = \{\langle x, F_A(x), T_A(x), I_A(x)\rangle : x \in X\},$$

which is characterized by a truth-membership function $T_A: X \to [0,1]$, an indeterminacy-membership function $I_A: X \to [0,1]$ and a falsity-membership function $F_A: X \to [0,1]$.

There is not restriction on the sum of $T_A(x)$, $I_A(x)$ and $F_A(x)$, so $0 \leq \sup T_A(x) \leq \sup I_A(x) \leq \sup F_A(x) \leq 3$.

Set- theoretic operations, for two single valued neutrosophic set,

$A_{NS} = \{<x, T_A(x), I_A(x), F_A(x)> | x \in X\}$ and $B_{NS} = \{<x, T_B(x), I_B(x), F_B(x)> | x \in X\}$ are given as;

1. The subset; $A_{NS} \subseteq B_{NS}$ if and only if

   $T_A(x) \leq T_B(x), I_A(x) \geq I_B(x), F_A(x) \geq F_B(x)$.

2. $A_{NS} = B_{NS}$ if and only if,

   $T_A(x) = T_B(x), I_A(x) = I_B(x), F_A(x) = F_B(x)$ for any $x \in X$.

3. The complement of $A_{NS}$ is denoted by $A_{NS}^o$ and is defined by

   $A_{NS}^o = \{<x, F_A(x), 1 - I_A(x), T_A(x)| x \in X\}$

4. The intersection

   $A_{NS} \cap B_{NS} = \{<x, \min\{T_A(x), T_B(x)\}, \max\{I_A(x), I_B(x)\}, \max\{F_A(x), F_B(x)\}>: x \in X\}$

5. The union

   $A_{NS} \cup B_{NS} = \{<x, \max\{T_A(x), T_B(x)\}, \min\{I_A(x), I_B(x)\}, \min\{F_A(x), F_B(x)\}>: x \in X\}$

A single valued neutrosophic number is denoted by $\tilde{A} = \langle T, I, F \rangle$ for convenience.

**Definition 2.2.** [15] Let $\tilde{A}_1 = \langle T_1, I_1, F_1 \rangle$ and $\tilde{A}_2 = \langle T_2, I_2, F_2 \rangle$ be two single valued neutrosophic number. Then, the operations for NNs are defined as below;

    i.    $\lambda \tilde{A} = \langle 1 - (1 - T_1)^\lambda, I_1^\lambda, F_1^\lambda \rangle$
    ii.   $\tilde{A}_1^\lambda = \langle T_1^\lambda, 1 - (1 - I_1)^\lambda, 1 - (1 - F_1)^\lambda \rangle$
    iii.  $\tilde{A}_1 + \tilde{A}_2 = \langle T_1 + T_2 - T_1 T_2, I_1 I_2, F_1 F_2 \rangle$

iv. $\tilde{A}_1 . \tilde{A}_2 = \langle T_1 T_2, I_1 + I_2 - I_1 I_2, F_1 + F_2 - F_1 F_2 \rangle$

where $\lambda > 0$.

**Definition 2.3.** [15] Let $\tilde{A}_1 = \langle T_1, I_1, F_1 \rangle$ be a single valued neutrosophic number. Then, the score function $s(\tilde{A}_1)$, accuracy function $a(\tilde{A}_1)$ and certainty function $c(\tilde{A}_1)$ of an SNN are defined as follows:

i. $s(\tilde{A}_1) = (T_1 + 1 - I_1 + 1 - F_1)/3$;
ii. $a(\tilde{A}_1) = T_1 - F_1$;
iii. $c(\tilde{A}_1) = T_1$

**Definition 2.4.** [15] Let $\tilde{A}_1 = \langle T_1, I_1, F_1 \rangle$ and $\tilde{A}_2 = \langle T_2, I_2, F_2 \rangle$ be two single valued neutrosophic number. The comparison method can be defined as follows:

i. if $s(\tilde{A}_1) > s(\tilde{A}_2)$, then $\tilde{A}_1$ is greater than $\tilde{A}_2$, that is, $\tilde{A}_1$ is superior to $\tilde{A}_2$, denoted by $\tilde{A}_1 > \tilde{A}_2$
ii. if $s(\tilde{A}_1) = s(\tilde{A}_2)$ and $a(\tilde{A}_1) > a(\tilde{A}_2)$, then $\tilde{A}_1$ is greater than $\tilde{A}_2$, that is, $\tilde{A}_1$ is superior to $\tilde{A}_2$, denoted by $\tilde{A}_1 < \tilde{A}_2$;
iii. if $s(\tilde{A}_1) = s(\tilde{A}_2)$, $a(\tilde{A}_1) = a(\tilde{A}_2)$ and $c(\tilde{A}_1) > c(\tilde{A}_2)$, then $\tilde{A}_1$ is greater than $\tilde{A}_2$, that is, $\tilde{A}_1$ is superior to $\tilde{A}_2$, denoted by $\tilde{A}_1 > \tilde{A}_2$;
iv. if $s(\tilde{A}_1) = s(\tilde{A}_2)$, $a(\tilde{A}_1) = a(\tilde{A}_2)$ and $c(\tilde{A}_1) = c(\tilde{A}_2)$, then $\tilde{A}_1$ is equal to $\tilde{A}_2$, that is, $\tilde{A}_1$ is indifferent to $\tilde{A}_2$, denoted by $\tilde{A}_1 = \tilde{A}_2$.

**Definition 2.4.** [6,14] Let X be a non-empty set. Then, a bipolar-valued fuzzy set, denoted by $A_{BF}$, is difined as;

$$A_{BF} = \{\langle x, \mu_B^+(x), \mu_B^-(x) \rangle : x \in X\}$$

where $\mu_B^+ : X \rightarrow [0,1]$ and $\mu_B^- : X \rightarrow [0,1]$. The positive membership degree $\mu_B^+(x)$ denotes the satisfaction degree of an element x to the property corresponding to $A_{BF}$ and the negative membership degree $\mu_B^-(x)$ denotes the satisfaction degree of x to some implicit counter property of $A_{BF}$.

## 3. Bipolar Neutrosophic Set

In this section, we introduce concept of bipolar neutrosophic set and its some operations including the score, certainty and accuracy functions to compare the bipolar neutrosophic sets. We also develop the bipolar neutrosophic weighted average operator ($A_w$) and bipolar neutrosophic weighted geometric operator ($G_w$) operator to aggregate the bipolar neutrosophic information. Some of it is quoted from [2,6,8,9,14,17,18,20,24,26].

**Definition 3.1.** A bipolar neutrosophic set $A$ in $X$ is defined as an object of the form

$$A = \{\langle x, T^+(x), I^+(x), F^+(x), T^-(x), I^-(x), F^-(x) \rangle : x \in X\},$$

where $T^+, I^+, F^+ : X \rightarrow [1,0]$ and $T^-, I^-, F^- : X \rightarrow [-1,0]$.

The positive membership degree $T^+(x), I^+(x), F^+(x)$ denotes the truth membership, indeterminate membership and false membership of an element $x \in X$ corresponding to a bipolar neutrosophic set $A$ and the negative membership degree $T^-(x), I^-(x), F^-(x)$ denotes the truth membership, indeterminate membership and false membership of an element $x \in X$ to some implicit counter-property corresponding to a bipolar neutrosophic set $A$.

**Example 3.2.** Let $X = \{x_1, x_2, x_3\}$. Then

$$A = \begin{Bmatrix} \langle x_1, 0.5, 0.3, 0.1, -0.6, -0.4, -0.01 \rangle, \\ \langle x_2, 0.3, 0.2, 0.7, -0.02, -0.003, -0.5 \rangle, \\ \langle x_3, 0.8, 0.05, 0.4, -0.1, -0.5, -0.06 \rangle \end{Bmatrix}$$

is a bipolar neutrosophic subset of $X$.

**Theorem 3.4.** A bipolar neutrosophic set is the generalization of a bipolar fuzzy set.
Proof: Suppose that $X$ is a bipolar neutrosophic set. Then by setting the positive components $I^+$, $F^+$ equals to zero as well as the negative components $I^-$, $F^-$ equals to zero reduces the bipolar neutrosophic set to bipolar fuzzy set.

**Definition 3.5.** Let $A_1 = \langle x, T_1^+(x), I_1^+(x), F_1^+(x), T_1^-(x), I_1^-(x), F_1^-(x) \rangle$ and
$A_2 = \langle x, T_2^+(x), I_2^+(x), F_2^+(x), T_2^-(x), I_2^-(x), F_2^-(x) \rangle$ be two bipolar neutrosophic sets. Then $A_1 \subseteq A_2$ if and only if

$$T_1^+(x) \leq T_2^+(x), \; I_1^+(x) \leq I_2^+(x), \; F_1^+(x) \geq F_2^+(x),$$

and

$$T_1^-(x) \geq T_2^-(x), \; I_1^-(x) \geq I_2^-(x), \; F_1^-(x) \leq F_2^-(x)$$

for all $x \in X$.

**Definition 3.6.** Let $A_1 = \langle x, T_1^+(x), I_1^+(x), F_1^+(x), T_1^-(x), I_1^-(x), F_1^-(x) \rangle$ and
$A_2 = \langle x, T_2^+(x), I_2^+(x), F_2^+(x), T_2^-(x), I_2^-(x), F_2^-(x) \rangle$ be two bipolar neutrosophic set. Then $A_1 = A_2$ if and only if

$$T_1^+(x) = T_2^+(x), \; I_1^+(x) = I_2^+(x), \; F_1^+(x) = F_2^+(x),$$

and

$$T_1^-(x) = T_2^-(x), \; I_1^-(x) = I_2^-(x), \; F_1^-(x) = F_2^-(x)$$

for all $x \in X$.

**Definition 3.7.** Let $A_1 = \langle x, T_1^+(x), I_1^+(x), F_1^+(x), T_1^-(x), I_1^-(x), F_1^-(x) \rangle$ and
$A_2 = \langle x, T_2^+(x), I_2^+(x), F_2^+(x), T_2^-(x), I_2^-(x), F_2^-(x) \rangle$ be two bipolar neutrosophic set. Then their union is defined as:

$$(A_1 \cup A_2)(x) = \left( \max(T_1^+(x), T_2^+(x)), \frac{I_1^+(x) + I_2^+(x)}{2}, \min((F_1^+(x), F_2^+(x)), \min(T_1^-(x), T_2^-(x)), \frac{I_1^-(x) + I_2^-(x)}{2}, \max((F_1^-(x), F_2^-(x))) \right)$$

for all $x \in X$.

**Example 3.8.** Let $X = \{x_1, x_2, x_3\}$. Then

$$A_1 = \begin{Bmatrix} \langle x_1, 0.5, 0.3, 0.1, -0.6, -0.4, -0.01 \rangle, \\ \langle x_2, 0.3, 0.2, 0.7, -0.02, -0.003, -0.5 \rangle, \\ \langle x_3, 0.8, 0.05, 0.4, -0.1, -0.5, -0.06 \rangle \end{Bmatrix}$$

and

$$A_2 = \begin{Bmatrix} \langle x_1, 0.4, 0.6, 0.3, -0.3, -0.5, -0.1 \rangle, \\ \langle x_2, 0.5, 0.1, 0.4, -0.2, -0.3, -0.3 \rangle, \\ \langle x_3, 0.2, 0.5, 0.6, -0.4, -0.6, -0.7 \rangle \end{Bmatrix}$$

are two bipolar neutrosophic sets in $X$.

Then their union is given as follows:

$$A_1 \cup A_2 = \begin{Bmatrix} \langle x_1, 0.5, 0.45, 0.1, -0.6, -0.5, -0.1 \rangle, \langle x_2, 0.5, 0.15, 0.7, -0.2, -0.1515, -0.5 \rangle, \\ \langle x_3, 0.8, 0.47, 0.6, -0.4, -0.55, -0.7 \rangle \end{Bmatrix}$$

**Definition 3.9.** Let $A_1 = \langle x, T_1^+(x), I_1^+(x), F_1^+(x), T_1^-(x), I_1^-(x), F_1^-(x) \rangle$ and $A_2 = \langle x, T_2^+(x), I_2^+(x), F_2^+(x), T_2^-(x), I_2^-(x), F_2^-(x) \rangle$ be two bipolar neutrosophic set. Then their intersection is defined as:

$$(A_1 \cap A_2)(x) = \left( \min(T_1^+(x), T_2^+(x)), \frac{I_1^+(x) + I_2^+(x)}{2}, \max((F_1^+(x), F_2^+(x)), \max(T_1^-(x), T_2^+(x)), \frac{I_1^-(x) + I_2^-(x)}{2}, \min((F_1^-(x), F_2^-(x)) \right)$$

for all $x \in X$.

**Definition 3.10.** Let $A = \left\{ \langle x, T^+(x), I^+(x), F^+(x), T^-(x), I^-(x), F^-(x) \rangle : x \in X \right\}$ be a bipolar neutrosophic set in $X$. Then the complement of $A$ is denoted by $A^c$ and is defined by

$$T_{A^c}^+(x) = \{1^+\} - T_A^+(x), \ I_{A^c}^+(x) = \{1^+\} - I_A^+(x), \ F_{A^c}^+(x) = \{1^+\} - F_A^+(x)$$

and

$$T_{A^c}^-(x) = \{1^-\} - T_A^-(x), \ I_{A^c}^-(x) = \{1^-\} - I_A^-(x), \ F_{A^c}^-(x) = \{1^-\} - F_A^-(x),$$

for all $x \in X$.

**Example 3.11.** Let $X = \{x_1, x_2, x_3\}$. Then

$$A = \begin{Bmatrix} \langle x_1, 0.5, 0.3, 0.1, -0.6, -0.4, -0.01 \rangle, \\ \langle x_2, 0.3, 0.2, 0.7, -0.02, -0.003, -0.5 \rangle, \\ \langle x_3, 0.8, 0.05, 0.4, -0.1, -0.5, -0.06 \rangle \end{Bmatrix}$$

is a bipolar neutrosophic set in $X$. Then the complement of $A$ is given as follows:

$$A^c = \begin{Bmatrix} \langle x_1, 0.5, 0.7, 0.9, -0.4, -0.6, -0.99 \rangle, \\ \langle x_2, 0.7, 0.8, 0.3, -0.08, -0.997, -0.5 \rangle, \\ \langle x_3, 0.2, 0.95, 0.6, -0.9, -0.5, -0.94 \rangle \end{Bmatrix}.$$

We will denote the set of all the bipolar neutrosophic sets (NBSs) in $X$ by $Q$. A bipolar neutrosophic number (NBN) is denoted by $\tilde{a} = \langle T^+, I^+, F^+, T^-, I^-, F^- \rangle$ for convenience.

**Definition 3.12.** Let $\tilde{a}_1 = \langle T_1^+, I_1^+, F_1^+, T_1^-, I_1^-, F_1^- \rangle$ and $\tilde{a}_2 = \langle T_2^+, I_2^+, F_2^+, T_2^-, I_2^-, F_2^- \rangle$ be two bipolar neutrosophic number. Then the operations for NNs are defined as below;

i. $\lambda\tilde{a}_1 = \langle 1 - (1 - T_1^+)^\lambda, (I_1^+)^\lambda, (F_1^+)^\lambda, -(-T_1^-)^\lambda, -(-I_1^-)^\lambda, -(1 - (1 - (-F_1^-))^\lambda) \rangle$
ii. $\tilde{a}_1^\lambda = \langle (T_1^+)^\lambda, 1 - (1 - I_1^+)^\lambda, 1 - (1 - F_1^+)^\lambda, -(1 - (1 - (-T_1^-))^\lambda), -(-I_1^-)^\lambda, -(-F_1^-)^\lambda \rangle$
iii. $\tilde{a}_1 + \tilde{a}_2 = \langle T_1^+ + T_2^+ - T_1^+ T_2^+, I_1^+ I_2^+, F_1^+ F_2^+, -T_1^- T_2^-, -(-I_1^- - I_2^- - I_1^- I_2^-), -(-F_1^- - F_2^- - F_1^- F_2^-) \rangle$
iv. $\tilde{a}_1 \cdot \tilde{a}_2 = \langle T_1^+ T_2^+, I_1^+ + I_2^+ - I_1^+ I_2^+, F_1^+ + F_2^+ - F_1^+ F_2^+, -(-T_1^- - T_2^- - T_2^- T_2^-), -I_1^- I_2^-, -F_1^- F_2^- \rangle$

where $\lambda > 0$.

**Definition 3.14.** Let $\tilde{a}_1 = \langle T_1^+, I_1^+, F_1^+, T_1^-, I_1^-, F_1^- \rangle$ be a bipolar neutrosophic number. Then, the score function $s(\tilde{a}_1)$, accuracy function $a(\tilde{a}_1)$ and certainty function $c(\tilde{a}_1)$ of an NBN are defined as follows:

i. $\tilde{s}(\tilde{a}_1) = (T_1^+ + 1 - I_1^+ + 1 - F_1^+ + 1 + T_1^- - I_1^- - F_1^-)/6$
ii. $\tilde{a}(\tilde{a}_1) = T_1^+ - F_1^+ + T_1^- - F_1^-$
iii. $\tilde{c}(\tilde{a}_1) = T_1^+ - F_1^-$

**Definition 3.15.** $\tilde{a}_1 = \langle T_1^+, I_1^+, F_1^+, T_1^-, I_1^-, F_1^- \rangle$ and $\tilde{a}_2 = \langle T_2^+, I_2^+, F_2^+, T_2^-, I_2^-, F_2^- \rangle$ be two bipolar neutrosophic number. The comparison method can be defined as follows:

i. if $\tilde{s}(\tilde{a}_1) > \tilde{s}(\tilde{a}_2)$, then $\tilde{a}_1$ is greater than $\tilde{a}_2$, that is, $\tilde{a}_1$ is superior to $\tilde{a}_2$, denoted by $a_1 > \tilde{a}_2$
ii. $\tilde{s}(\tilde{a}_1) = \tilde{s}(\tilde{a}_2)$ and $\tilde{a}(\tilde{a}_1) > \tilde{a}(\tilde{a}_2)$, then $\tilde{a}_1$ is greater than $\tilde{a}_2$, that is, $\tilde{a}_1$ is superior to $\tilde{a}_2$, denoted by $\tilde{a}_1 < \tilde{a}_2$;
iii. if $\tilde{s}(\tilde{a}_1) = \tilde{s}(\tilde{a}_2)$, $\tilde{a}(\tilde{a}_1) = \tilde{a}(\tilde{a}_1)$ and $\tilde{c}(\tilde{a}_1) > \tilde{c}(\tilde{a}_2)$, then $\tilde{a}_1$ is greater than $\tilde{a}_2$, that is, $\tilde{a}_1$ is superior to $\tilde{a}_2$, denoted by $\tilde{a}_1 > \tilde{a}_2$;
iv. if $\tilde{s}(\tilde{a}_1) = \tilde{s}(\tilde{a}_2)$, $\tilde{a}(\tilde{a}_1) = \tilde{a}(\tilde{a}_2))$ and $\tilde{c}(\tilde{a}_1) = \tilde{c}(\tilde{a}_2)$, then $\tilde{a}_1$ is equal to $\tilde{a}_2$, that is, $\tilde{a}_1$ is indifferent to $\tilde{a}_2$, denoted by $\tilde{a}_1 = \tilde{a}_2$.

Based on the study given in [15,20] we define some weighted aggregation operators related to bipolar neutrosophic sets as follows;

**Definition 3.16.** Let $\tilde{a}_j = \langle T_j^+, I_j^+, F_j^+, T_j^-, I_j^-, F_j^- \rangle (j = 1,2,...,n)$ be a family of bipolar neutrosophic numbers. A mapping $A_\omega: Q_n \to Q$ *is called* bipolar neutrosophic weighted average operator if it satisfies

$A_w(\tilde{a}_1, \tilde{a}_2, ..., \tilde{a}_n) = \sum_{j=1}^n \omega_j \tilde{a}_j$

$= \langle 1 - \prod_{j=1}^n (1 - T_j^+)^{\omega_j}, \prod_{j=1}^n I_j^{+\omega_j}, \prod_{j=1}^n F_j^{+\omega_j}, -\prod_{j=1}^n (-T_j^-)^{\omega_j}, -(1 - \prod_{j=1}^n (1 - (-I_j^-))^{\omega_j}), -(1 - \prod_{j=1}^n (1 - (-F_j^-))^{\omega_j}) \rangle$

where $\omega_j$ is the weight of $\tilde{a}_j$ ($j = 1,2, \ldots, n$), $\omega_j \in [0,1]$ and $\sum_{j=1}^{n} \omega_j = 1$.

Based on the study given in [15,20] we give the theorem related to bipolar neutrosophic sets as follows;

**Theorem 3.17.** Let $\tilde{a}_j = \langle T^+_j, I^+_j, F^+_j, T^-_j, I^-_j, F^-_j \rangle$ ($j = 1,2, \ldots, n$) be a family of bipolar neutrosophic numbers. Then,

i. If $\tilde{a}_j = \tilde{a}$ for all $j = 1,2, \ldots, n$ then, $A_w(\tilde{a}_1, \tilde{a}_2, \ldots, \tilde{a}_n) = \tilde{a}$
ii. $\min_{j=1,2,\ldots,n}\{\tilde{a}_j\} \leq A_w(\tilde{a}_1, \tilde{a}_2, \ldots, \tilde{a}_n) \leq \max_{j=1,2,\ldots,n}\{\tilde{a}_j\}$
iii. If $\tilde{a}_j \leq \tilde{a}_j^*$ for all $j = 1,2, \ldots, n$ then, $A_w(\tilde{a}_1, \tilde{a}_2, \ldots, \tilde{a}_n) \leq A_w(\tilde{a}_1^*, \tilde{a}_2^*, \ldots, \tilde{a}_n^*)$

Based on the study given in [15,20] we define some weighted aggregation operators related to bipolar neutrosophic sets as follows;

**Definition 3.18.** Let $\tilde{a}_j = \langle T^+_j, I^+_j, F^+_j, T^-_j, I^-_j, F^-_j \rangle$ ($j = 1,2, \ldots, n$) be a family of bipolar neutrosophic numbers. A mapping $G_\omega: \mathcal{Q}_n \to \mathcal{Q}$ is called bipolar neutrosophic weighted geometric operator if it satisfies

$$G_w(\tilde{a}_1, \tilde{a}_2, \ldots, \tilde{a}_n) = \prod_{j=1}^{n} \tilde{a}_j^{\omega_j} = \langle \prod_{j=1}^{n} T_j^{+\omega_j}, 1 - \prod_{j=1}^{n}(1 - I_j^+)^{\omega_j}, 1 - \prod_{j=1}^{n}(1 - F_j^+)^{\omega_j}, -(1 - \prod_{j=1}^{n}(1 - (-T_j^-))^{\omega_j}), -\prod_{j=1}^{n}(-I_j^-)^{\omega_j}, -\prod_{j=1}^{n}(-F_j^{-\omega_j}) \rangle$$

where $\omega_j$ is the weight of $\tilde{a}_j$ ($j = 1,2, \ldots, n$), $\omega_j \in [0,1]$ and $\sum_{j=1}^{n} \omega_j = 1$.

Based on the study given in [15,20] we give the theorem related to bipolar neutrosophic sets as follows;

**Theorem 3.19.** Let $\tilde{a}_j = \langle T^+_j, I^+_j, F^+_j, T^-_j, I^-_j, F^-_j \rangle$ ($j = 1,2, \ldots, n$) be a family of bipolar neutrosophic numbers. Then,

i. If $\tilde{a}_j = \tilde{a}$ for all $j = 1,2, \ldots, n$ then, $G_w(\tilde{a}_1, \tilde{a}_2, \ldots, \tilde{a}_n) = \tilde{a}$
ii. $\min_{j=1,2,\ldots,n}\{\tilde{a}_j\} \leq G_w(\tilde{a}_1, \tilde{a}_2, \ldots, \tilde{a}_n) \leq \max_{j=1,2,\ldots,n}\{\tilde{a}_j\}$
iii. If $\tilde{a}_j \leq \tilde{a}_j^*$ for all $j = 1,2, \ldots, n$ then, $G_w(\tilde{a}_1, \tilde{a}_2, \ldots, \tilde{a}_n) \leq G_w(\tilde{a}_1^*, \tilde{a}_2^*, \ldots, \tilde{a}_n^*)$

Note that the aggregation results are still NBNs

## 4. NBN- Decision Making Method

In this section, we develop an approach based on the $A_w$ (or $G_w$) operator and the above ranking method to deal with multiple criteria decision making problems with bipolar neutrosophic information.

Suppose that $A = \{A_1, A_2, \ldots, A_m\}$ and $C = \{C_1, C_2, \ldots, C_n\}$ is the set of alternatives and criterions or attributes, respectively. Let $\omega = (\omega_1, \omega_2, \ldots, \omega_n)^T$ be the weight vector of attributes, such that $\sum_{j=1}^{n} \omega_j = 1$, $\omega_j \geq 0$ $(j = 1,2, \ldots, n)$ and $\omega_j$ refers to the weight of attribute $C_j$. An alternative on criterions is evaluated by the decision maker, and the evaluation values are represented by the form of bipolar neutrosophic numbers. Assume that $(\tilde{a}_{ij})_{m \times n} = \left(\langle T_{ij}^+, I_{ij}^+, F_{ij}^+, T_{ij}^-, I_{ij}^-, F_{ij}^- \rangle\right)_{m \times n}$ is the decision matrix provided by the decision maker; $\tilde{a}_{ij}$ is a bipolar neutrosophic number for alternative $A_i$ associated with the criterions $C_j$. We have the conditions $T_{ij}^+, I_{ij}^+, F_{ij}^+, T_{ij}^-, I_{ij}^-$ and $F_{ij}^- \in [0,1]$ such that $0 \leq T_{ij}^+ + I_{ij}^+ + F_{ij}^+ - T_{ij}^- - I_{ij}^- - F_{ij}^- \leq 6$ for $i = 1,2, \ldots, m$ and $j = 1,2, \ldots, n$.

Now, we can develop an algorithm as follows;

**Algorithm**

**Step1.** Construct the decision matrix provided by the decision maker as;
$$(\tilde{a}_{ij})_{m \times n} = \left(\langle T_{ij}^+, I_{ij}^+, F_{ij}^+, T_{ij}^-, I_{ij}^-, F_{ij}^- \rangle\right)_{m \times n}$$

**Step 2.** Compute $\tilde{a}_i = A_w(\tilde{a}_{i1}, \tilde{a}_{i2}, \ldots, \tilde{a}_{in})$ (or $G_w(\tilde{a}_{i1}, \tilde{a}_{i2}, \ldots, \tilde{a}_{in})$) for each $i = 1,2, \ldots, m$.

**Step 3.** Calculate the score values of $\tilde{s}(\tilde{a}_1)$ $(i = 1,2, \ldots, m.)$ for the collective overall bipolar neutrosophic number of $\tilde{a}_i$ $(i = 1,2, \ldots, m.)$

**Step 4.** Rank all the software systems of $\tilde{a}_i$ $(i = 1,2, \ldots, m.)$ according to the score values

Now, we give a numerical example as follows;

**Example 4.1.** Let us consider decision making problem adapted from Xu and Cia [20]. A customer who intends to buy a car. Four types of cars (alternatives) $A_i (i = 1,2,3,4)$ are available. The customer takes into account four attributes to evaluate the alternatives; $C_1 =$ Fuel economy; $C_2 =$ Aerod; $C_3 =$ Comfort; $C_4 =$ Safety and use the bipolar neutrosophic values to evaluate the four possible alternatives $A_i (i = 1, 2, 3, 4)$ under the above four attributes. Also, the weight vector of the attributes $C_j (j = 1,2,3,4)$ is $\omega = (\frac{1}{2}, \frac{1}{4}, \frac{1}{8}, \frac{1}{8})^T$. Then,

**Algorithm**

**Step 1.** Construct the decision matrix provided by the customer as;

Table 1: Decision matrix given by customer

|   | $C_1$ | $C_2$ | $C_3$ | $C_4$ |
|---|---|---|---|---|
| $A_1$ | ⟨0.5,0.7,0.2,−0.7,−0.3,−0.6⟩ | ⟨0.4,0.4,0.5,−0.7,−0.8,−0.4⟩ | ⟨0.7,0.7,0.5,−0.8,−0.7,−0.6⟩ | ⟨0.1,0.5,0.7,−0.5,−0.2,−0.8⟩ |
| $A_2$ | ⟨0.9,0.7,0.5,−0.7,−0.7,−0.1⟩ | ⟨0.7,0.6,0.8,−0.7,−0.5,−0.1⟩ | ⟨0.9,0.4,0.6,−0.1,−0.7,−0.5⟩ | ⟨0.5,0.2,0.7,−0.5,−0.1,−0.9⟩ |
| $A_3$ | ⟨0.3,0.4,0.2,−0.6,−0.3,−0.7⟩ | ⟨0.2,0.2,0.2,−0.4,−0.7,−0.4⟩ | ⟨0.9,0.5,0.5,−0.6,−0.5,−0.2⟩ | ⟨0.7,0.5,0.3,−0.4,−0.2,−0.2⟩ |
| $A_4$ | ⟨0.9,0.7,0.2,−0.8,−0.6,−0.1⟩ | ⟨0.3,0.5,0.2,−0.5,−0.5,−0.2⟩ | ⟨0.5,0.4,0.5,−0.1,−0.7,−0.2⟩ | ⟨0.4,0.2,0.8,−0.5,−0.5,−0.6⟩ |

**Step 2.** Compute $\tilde{a}_i = A_w(\tilde{a}_{i1}, \tilde{a}_{i2}, \tilde{a}_{i3}, \tilde{a}_{i4})$ for each $i = 1,2,3,4$ as;

$$\tilde{a}_1 = \langle 0.471, 0.583, 0.329, -0.682, -0.531, -0.594 \rangle$$

$$\tilde{a}_2 = \langle 0.839, 0.536, 0.600, -0.526, -0.608, -0.364 \rangle$$

$$\tilde{a}_3 = \langle 0.489, 0.355, 0.235, -0.515, -0.447, -0.544 \rangle$$

$$\tilde{a}_4 = \langle 0.751, 0.513, 0.266, -0.517, -0.580, -0.221 \rangle$$

**Step 3.** Calculate the score values of $\tilde{s}(\tilde{a}_1)$ ($i = 1,2,3,4$) for the collective overall bipolar neutrosophic number of $\tilde{a}_i$ ($i = 1,2,...,m.$) as;

$$\tilde{s}(\tilde{a}_1) = 0.50$$

$$\tilde{s}(\tilde{a}_2) = 0.52$$

$$\tilde{s}(\tilde{a}_3) = 0.56$$

$$\tilde{s}(\tilde{a}_4) = 0.54$$

**Step 4.** Rank all the software systems of $A_i$ ($i = 1,2,3,4.$) according to the score values as;

$$A_3 \succ A_4 \succ A_2 \succ A_1$$

and thus $A_3$ is the most desirable alternative.

## 5. Conclusions

This paper presented a bipolar neutrosophic set and its score, certainty and accuracy functions. Then, the $A_w$ and $G_w$ operators were proposed to aggregate the bipolar neutrosophic information. Furthermore, based on the $A_w$ and $G_w$ operators and the score, certainty and accuracy functions, we have developed a bipolar neutrosophic multiple criteria decision-making approach, in which the evaluation values of alternatives on the attributes take the form of bipolar neutrosophic numbers. The $A_w$ and $G_w$ operators are utilized to aggregate the bipolar neutrosophic information corresponding to each alternative to obtain the collective overall values of the alternatives, and then the alternatives are ranked according to the values of the score, certainty and accuracy functions to select the most desirable one(s). Finally, a numerical example of the method was given to demonstrate the application and effectiveness of the developed method.

**References**

[1]  K. Atanassov, Intuitionistic fuzzy sets. Fuzzy Sets and Systems, 20  (1986) 87–96.
[2] P. Bosc, O. Pivert, On a fuzzy bipolar relational algebra, Information Sciences, 219 (2013) 1–16.
[3] J. Chen, S. Li, S. Ma, and X. Wang, $m$-Polar Fuzzy Sets: An Extension of Bipolar Fuzzy Sets, The Scientific World Journal, (2014) http://dx.doi.org/10.1155/2014/416530.
[4] H. D. Cheng and Y. Guo, A new neutrosophic approach to image thresholding, New Mathematics and Natural Computation, 4(3) (2008) 291–308.
[5] Y. Guo and H. D. Cheng, New Neutrosophic Approach to Image Segmentation, Pattern Recognition, 42, (2009), 587–595.
[6] A. Kharal, A neutrosophic multicriteria decision making method, New Mathematics & Natural Computation, 2013.
[7] M. K. Kang and J. G. Kang, Bipolar fuzzy set theory applied to sub-semigroups with operators in semigroups. J. Korean Soc. Math. Educ. Ser. B Pure Appl. Math., 19/1 (2012) 23-35.
[8] K. M. Lee, Bipolar-valued fuzzy sets and their operations. Proc. Int. Conf. on Intelligent Technologies, Bangkok, Thailand  (2000) 307-312.
[9] K. J. Lee, Bipolar fuzzy subalgebras and bipolar fuzzy ideals of BCK/BCI-algebras,  Bull. Malays. Math. Sci. Soc., 32/3 (2009) 361-373.
[10] P. Liu and Y. Wang, Multiple attribute decision-making method based on single-valued neutrosophic normalized weighted Bonferroni mean, Neural Computing and Applications, 25 7/8 (2014)  2001-2010.
[11] P. Liu L. Shi, The generalized hybrid weighted average operator based on interval neutrosophic hesitant set and its application to multiple attribute decision making, Neural Computing and Applications, 26 /2 (2015) 457-471.
[12] P. Majumdar and S.K. Samanta, On similarity and entropy of neutrosophic sets, J. Intell. Fuzzy Syst. 26/3 (2014) 1245–1252.
[13]  S.K.Majumder, Bipolar Valued Fuzzy Sets in Γ-Semigroups, Mathematica Aeterna, 2/3 (2012) 203 – 213.
[14]  S.V. Manemaran B. Chellappa, Structures on Bipolar Fuzzy Groups and Bipolar Fuzzy D-Ideals under (T, S) Norms, International Journal of Computer Applications, 9/12, 7-10.
[15] J.J. Peng, J.Q. Wang, J. Wang, H.Y. Zhang and X.H. Chen, Simplified neutrosophic sets and their applications in multi-criteria group decision-making problems, Int. J. Syst. Sci. (2015) DOI:10.1080/00207721.2014.994050.
[16] R. Sahin and A. Kucuk, Subsethood measure for single valued neutrosophic sets, Journal of Intelligent and Fuzzy Systems, (2014)  DOI: 10.3233/IFS-141304.
[17] F. Smarandache, A Unifying Field in Logics. Neutrosophy : Neutrosophic Probability, Set and Logic, Rehoboth: American Research Press,1999.
[18] H. Wang, F. Smarandache, Y.Q. Zhang and R. Sunderraman Single valued neutrosophic sets, Multispace and Multistructure, 4  (2010)  10-413.
[19] H. Wang, F. Smarandache, Y.Q. Zhang and R. Sunderraman, Interval neutrosophic sets and logic: theory and applications in computing, (2005)  Hexis, Arizona.

[20] Z. Xu, X. Cai, Intuitionistic Fuzzy Information Aggregation Theory and Applications, Springer, Science Press, Heidelberg New York Dordrecht London, 2011.

[21] J. Ye Multicriteria decision-making method using the correlation coefficient under single-valued neutrosophic environment, International Journal of General Systems 42/4 (2013) 386-394.

[22] J. Ye, Similarity measures between interval neutrosophic sets and their applications in Multi-criteria decision-making. Journal of Intelligent and Fuzzy Systems, 26 (2014) 165-172.

[23] J. Ye, Single valued neutrosophic cross-entropy for multi-criteria decision making problems, Appl. Math. Model. 38 /3 (2014) 1170–1175.

[24] J. Ye, Trapezoidal neutrosophic set and its application to multiple attribute decision making, Neural Computing and Applications, (2014) DOI: 10.1007/s00521-014-1787-6

[25] J. Ye, Some aggregation operators of interval neutrosophic linguistic numbers for multiple attribute decision making, Journal of Intelligent & Fuzzy Systems 27 (2014) 2231-2241

[26] L.A. Zadeh, Fuzzy sets, Inf. Control, 8 (1965) 338–353.

[27] H.Y. Zhang, J.Q. Wang, X.H. Chen, Interval neutrosophic sets and their application in multicriteria decision making problems. The Scientific World Journal. (2014) DOI: 10.1155/2014/645953.

[28] M. Zhang, L. Zhang, and H.D. Cheng. A Neutrosophic Approach to Image Segmentation based on Watershed Method, Signal Processing 5/ 90 (2010) 1510-1517.

[29] M.Zhou, S. LI, Application of Bipolar Fuzzy Sets in Semirings, Journal of Mathematical Research with Applications, Vol. 34/ 1 (2014) 61-72